\documentstyle[amssymb]{amsart}
\evensidemargin=\oddsidemargin
\def\Z{{\bold Z}}
\def\R{{\bold R}}
\def\bC{{\bold C}}
\def\bH{{\bold H}}
\def\bF{{\bold F}}
\def\a{{\frak{a}}}
\def\f{{\frak{f}}}
\def\g{{\frak{g}}}
\def\h{{\frak{h}}}
\def\fk{{\frak{k}}}

\def\m{{\frak{m}}}
\def\n{{\frak{n}}}
\def\gl{{\frak{gl}}}

\def\s{{\frak{s}}}
\def\so{{\frak{so}}}
\def\su{{\frak{su}}}
\def\sl{{\frak{sl}}}
\def\fu{{\frak{u}}}

\def\v{{\frak{v}}}
\def\z{{\frak{z}}}
\def\fsp{{\frak{sp}}}
\def\gr{{\g^\R}}
\def\rg{{\g_0}}
\def\S{{\Sigma}}
\def\la{{\langle}}
\def\ip{{\la\,,\,\rangle}}
\def\al{{\alpha}}
\def\b{{\beta}}
\def\bb{{\bar{\b}}}

\def\half{{\tfrac 12}}
\def\fourth{{\tfrac 14}}
\def\H{{\bold H}}
\def\PP{{P\oplus P}}
\def\i{{\bold i}}
\def\j{{\bold j}}
\def\k{{\bold k}}
\def\ep{{\epsilon}}
\def\om{{\omega}} 
\def\ra{{\rightarrow}}
\def\PP{{P\oplus P}}
\newcommand\tr{\operatorname{tr}}
\newcommand\SL{\operatorname{SL}}
\newcommand\SU{\operatorname{SU}}
\newcommand\U{\operatorname{U}}
\newcommand\SO{\operatorname{SO}}
\renewcommand\Sp{\operatorname{Sp}}
\renewcommand\O{\operatorname{O}}

\newcommand\Id{\operatorname{Id}}
\newcommand\ad{\operatorname{ad}}
\newcommand\ric{\operatorname{Ric}}
\newcommand\image{\operatorname{image}}
\theoremstyle{plain}
\newtheorem{theorem}{Theorem}
\newtheorem{lemma}[theorem]{Lemma}
\newtheorem{proposition}[theorem]{Proposition}
\newtheorem{corollary}[theorem]{Corollary}

\theoremstyle{definition}
\newtheorem{definition}[theorem]{Definition}
\newtheorem{notation}[theorem]{Notation}
\newtheorem{example}[theorem]{Example}

\theoremstyle{remark}
\newtheorem{remark}[theorem]{Remark}

\numberwithin{equation}{section}
\numberwithin{theorem}{section}
\begin{document}
\title[Homogeneous Einstein Metrics]
{New Homogeneous Einstein Metrics \\ of Negative Ricci Curvature}
\author[C. S. Gordon]{Carolyn S. Gordon}
\address{Department of Mathematics, Dartmouth College, Hanover, NH 03755}
\email{Carolyn.S.Gordon@@Dartmouth.edu}
\thanks{1991 {\it Mathematical Subject Classification.} 
Primary 53C30.}
\author[M. M. Kerr]{Megan M. Kerr}
\address{Department of Mathematics, Wellesley College, 106 Central St., 
Wellesley, MA 02481}
\email{mkerr@@wellesley.edu}
\date{\today}

\begin{abstract} We construct new homogeneous Einstein spaces with negative
Ricci curvature in two ways: First, we give a method for classifying and constructing
a class of rank one Einstein solvmanifolds whose derived algebras are two-step nilpotent.  As an
application, we describe an explicit continuous family of ten-dimensional Einstein
manifolds with a two-dimensional parameter space, including a continuous subfamily of manifolds with negative
sectional curvature. Secondly, we obtain new examples of  non-symmetric Einstein solvmanifolds
by modifying the algebraic structure of non-compact irreducible symmetric spaces of rank
greater than one, preserving the (constant) Ricci curvature.  
\end{abstract}
\maketitle

\section*{Introduction}
    
The primary goal of this paper is the construction of new examples of
homogeneous Einstein manifolds of negative Ricci curvature.  The classical examples of
Einstein metrics of negative Ricci curvature are the symmetric spaces of non-compact
type.  A number of other examples are known, e.g. \cite{Al}, \cite{D}, \cite{GP-SV},
\cite{P-S}, and \cite{Wt2} (see \cite{Hb2} for more extensive references).  

A Riemannian solvmanifold is a solvable Lie group $S$ together with a
left-invariant Riemannian metric.  All known examples of homogeneous Einstein manifolds of
negative Ricci curvature are Riemannian solvmanifolds of Iwasawa type.   See Section 1
for a precise definition; here we just note that a solvable Lie group of Iwasawa type is
a semi-direct product of an abelian group $A$ with a nilpotent normal subgroup $N$, the
nilradical.  We restrict our attention in this paper to solvmanifolds of Iwasawa type.  

After providing the necessary background in Section 1, we consider in Section 2 the class
of solvmanifolds for which $N$ is two-step nilpotent, $A$ is one-dimensional, and the
action of $A$ on $N$ is of a particular type.  The solvmanifolds in this class are sometimes
referred to as Carnot solvmanifolds. This class of solvmanifolds includes (but is not limited
to) the harmonic manifolds of negative Ricci curvature constructed by Damek and Ricci
\cite{DR}. P. Eberlein and J. Heber \cite{EH2} gave a method for classifying all Einstein
solvmanifolds of this type; see Theorem \ref{ec} below.  For completeness, we include the
proof of this result (which we observed independently).  We then give an explicit
classification in low dimensions with several new examples, including, in Section 3, an
explicit continuous family of ten-dimensional Einstein manifolds with a two-dimensional
parameter space. We prove that within the deformation there is a continuous family of
non-isometric Einstein solvmanifolds with negative sectional curvature.

In Section 4, we obtain new examples of non-symmetric Einstein solvmanifolds by 
modifying the algebraic structure of the non-compact irreducible symmetric spaces of rank
greater than one.  

The study of homogeneous Einstein metrics is currently very active.  While this work was
in progress, we learned of recent work on Einstein solvmanifolds contained in
the habilitation of Jens Heber \cite{Hb2} and the
thesis of Denise Hengesch \cite{Hn}.  Heber's extensive habilitation
is a ground-breaking work, which allows Einstein solvmanifolds (which are  non-compact) to be
viewed as critical points of a certain functional and which delves deeply into questions of
existence and uniqueness.  Although our results in Section 4 were obtained independently of
Heber's habilitation, they share a common theme with his work:  Beginning with a given
Einstein solvmanifold, we modify not the metric but rather the bracket structure on the
associated Lie algebra in order to construct new examples of Einstein solvmanifolds.  Our
work in Sections 2 and 3 was motivated directly by Heber's habilitation.  In his study of the
moduli space of Einstein solvmanifolds, he shows that the moduli space of Einstein
solvmanifolds near a given one may have very large dimension.  In some cases one can compute
the dimension explicitly; in other cases one gets lower bounds on the dimension.  However,
his computation of the dimension does not give explicit examples but rather guarantees that
if there is one Einstein solvmanifold of a particular type, then (depending on the dimension
and structure) one  may be guaranteed a large moduli space. Our classification results and
explicit examples in Sections 2 and 3 complement Heber's work.  In particular, one can show
using Heber's work together with Theorem \ref{ec} that the lowest dimension in
which continuous families of Einstein Carnot solvmanifolds can exist is dimension 10;
i.e., our examples in Section 3 have the lowest possible dimension.

We wish to thank Wolfgang Ziller for helpful and stimulating discussions.  We also
thank Patrick Eberlein for  bringing to our attention the thesis of Sven Leukert \cite{L}
which was helpful in showing that the deformation in Section 3 contains negatively curved
examples.

%
%
%
\section{Preliminaries}
A solvable Lie group $S$ together with a left-invariant Riemannian metric $g$ is called a
{\it Riemannian solvmanifold}.  The left-invariant Riemannian metric $g$ on $S$ defines
an inner product $\langle \,,\,\rangle$ on the Lie algebra $\s$ of $S$, and conversely,
any inner product on $\s$ gives rise to a unique left-invariant metric on $S$.  

We will refer to a Lie algebra endowed with an inner product as a {\it metric Lie algebra}.
Thus we have a one-to-one correspondence between simply-connected Riemannian solvmanifolds and
solvable metric Lie algebras.  We will say two metric Lie algebras $(\s,\ip)$ and $(\s ',\ip
')$ are {\it isomorphic} if there exists a linear map $\tau :\s\to\s '$ which is both an
inner product space isometry and a Lie algebra isomorphism.  

\begin{notation}\label{iwasawa} Let  $(S,g)$ be a simply-connected Riemannian
solvmanifold, and let $\s$ be the associated solvable metric Lie algebra.  Write $\s=\a
\oplus\n$ with $\n$ the nilradical of $\s$. The metric Lie algebra $(\s,\ip)$ and the 
solvmanifold $(S,g)$ are said to be of {\it Iwasawa type} if the following conditions
hold:
\begin{itemize}
\item[(i)] $\a$ is abelian;
\item[(ii)] All $\ad(A)$, $A\in\a$, are symmetric relative to the inner product $\ip$ on
$\s$ and non-zero when $A\neq 0$;
\item[(iii)] For some $A\in\a$, the restriction of $\ad(A)$ to $\n$ is positive-definite.
\end{itemize}
When (i) holds, the dimension of $\a$ is called the {\it algebraic rank} of $\s$ and of
$S$.
\end{notation}

The following condition for two Riemannian solvmanifolds of Iwasawa type to be isometric
is a special case of Theorem 5.2 in \cite{GW}.

\begin{proposition}\label{iso}  Two simply-connected solvmanifolds of Iwasawa type are
isometric if and only if the associated metric Lie algebras are isomorphic.
\end{proposition}

 We next recall the formula for the Ricci tensor of the solvmanifold $(S,g)$, viewed as a
bilinear form on $\s$; see Besse \cite{Bes} for details.  

\begin{notation}\label{U}  Let $\{X_1,\dots,X_n\}$ be an orthonormal basis for $\s$.  
Define $U: \s \times \s \to \s$ by 
$$\langle U(X,Y),Z\rangle =\half\langle [Z,X],Y\rangle + \half\langle [Z,Y],X\rangle $$
for $X$, $Y$, $Z\in\s$.  Set
$$H = \sum_i U(X_i,X_i).$$
Observe that $H$ must lie in $\a$ and that 
$$\langle H,X\rangle = \tr\,\ad (X)$$
for all $X\in\s$.  As an aside, we note that $H$ is the mean curvature vector field for
the embedding of the nilradical $N$ in $S$. Thus, we will refer to $H$ as the mean
curvature vector field.
\end{notation}

Let $B$ denote the Killing form of $\s$.  For $X,Y\in\s$, the Ricci curvature $\ric(X,Y)$
is given by 
\begin{equation}\label{ricci}\ric(X,Y) = -\tfrac 12 \sum_i \langle
[X,X_i],[Y,X_i]\rangle -
\tfrac 12 B(X,Y) + \tfrac 14 \sum_{i,j}
\langle[X_i,X_j],X\rangle\langle[X_i,X_j],Y\rangle -
\langle U(X,Y),H\rangle.
\end{equation}
In particular,
\begin{equation}\label{ricci2}\ric(X,X) =-\tfrac 12 \tr\,\ad(X)\circ\ad(X)^* -
\tfrac 12 B(X,X) + \tfrac 14 \sum_{i,j} \langle [X_i,X_j], X\rangle ^2 -\langle
[H,X],X\rangle.
\end{equation}

As pointed out in \cite{Wt2} and easily seen from equation (\ref{ricci}), if the manifold
$(S,g)$ is Einstein, then so is the solvmanifold associated with the subalgebra
$\s_0 =\R H +\n$ of $\s$ with the induced inner product.  Conversely, Heber \cite{Hb2}
showed how to construct Einstein solvmanifolds of higher algebraic rank from a given
Einstein solvmanifold of Iwasawa type and of rank one.  Thus in Section 2, we will
focus on a class  of solvmanifolds  of algebraic rank one.  In Section 4, it will be more
convenient to work with higher rank solvmanifolds, however.

\begin{proposition}\label{eig} \cite{Hb2}  Let $(S,g)$ be an Einstein solvmanifold of
Iwasawa type and let $H\in\s$ be the mean curvature vector field defined in Notation
\ref{U}. Then for some multiple $\lambda H$ of $H$, the eigenvalues of $\ad(\lambda
H)_{|\n}$ are positive integers.  Moreover, the distinct eigenvalues have no common
divisors.
\end{proposition}

If $\mu_1 < \mu_2 <\dots <\mu_m$ are the distinct eigenvalues of $\ad(\lambda H)_{|\n}$
with multiplicities $d_1,\dots,d_m$, respectively, Heber refers to the data $(\mu_1,\dots,\mu_m
;d_1,\dots,d_m)$ as the {\it eigenvalue type} of the Einstein solvmanifold.  In Section 
\ref{twostep} we will consider Einstein solvmanifolds of Iwasawa type and of algebraic
rank one  with eigenvalue type $(1,2; d_1,d_2)$ where $d_1$ and $d_2$ are arbitrary. 
This is a rich class of Riemannian solvmanifolds which contains the negatively curved
rank one symmetric spaces and all the  non-compact harmonic manifolds constructed by
Damek and Ricci \cite{DR}. In Section 4 we will consider Einstein solvmanifolds of
the same eigenvalue types as the higher rank symmetric spaces.

\section{Two-Step Examples}\label{twostep}

In the notation following Proposition \ref{eig}, we now consider Einstein solvmanifolds
of algebraic rank one and eigenvalue type $(1,2;r,s)$ where the multiplicities $r$ and
$s$ of the eigenvalues are arbitrary.  All solvmanifolds considered here will be of
Iwasawa type.  We use the terminology
of Notations \ref{iwasawa} and \ref{U}.  Letting $\v$ and $\z$ be the eigenspaces of
$\ad(\lambda H)_{|\n}$ corresponding to the eigenvalues 1 and 2, respectively, the Jacobi
identity implies that $[\v,\v]\subset\z$ and that $\z$ is central in $\n$.  Thus $\n$ is
either  two-step nilpotent or abelian.  (As we will see in Theorem \ref{ec}, $\n$ must in fact be
two-step nilpotent, not abelian, in order for the Einstein condition to hold.)  Since
$(\s,\ip)$ is of Iwasawa type, $\ad(H)$ is symmetric and thus $\v$ and $\z$ are
orthogonal.

\begin{notation}\label{vz} (i) Let $(\n, \ip)$ be a two-step nilpotent metric Lie algebra, let
$\z$ be a central subspace of $\n$ containing the derived algebra $[\n,\n]$, and let $\v =
\z^{\perp}$ relative to $\langle\cdot\,,\,\cdot\rangle$.  We can then define a linear map
$j : \z \ra \so(\v, \langle\cdot\,,\,\cdot\rangle)$ by
\begin{equation}\label{j}
\langle j(Z) X, Y \rangle = \langle [X, Y], Z \rangle \text{ for }\,  
X, \, Y \in \v, \, Z \in \z. 
\end{equation}
Equivalently,
\begin{equation}\label{j2} j(Z)X=\ad (X)^*(Z). \end{equation}

(ii) Conversely, given any two finite dimensional real inner product spaces
${\frak v}$ and ${\frak z}$ along with a linear map $j : {\frak z} \ra \so({\frak v})$,
we can define a metric Lie alegbra $\n$ as the inner product space direct sum of ${\frak
v}$ and ${\frak z}$ with the alternating bilinear bracket map  $[\cdot, \cdot] : \n
\times \n \ra 
\z$ defined by declaring $\z$ to be central in $\n$ and using (\ref{j}) to define $[X, Y]$  for
$X, Y \in {\frak v}$. Then $\n$ is two-step nilpotent if $j$ is non-zero.

(iii) Given the data $(\v,\z,j)$ and thus a two-step nilpotent (or abelian if $j=0$)
metric Lie algebra $(\n, \ip)$ as in (ii), we can further define a Riemannian
solvmanifold $(S,g)$ of algebraic rank one and eigenvalue type $(1,2;r,s)$ where
$r=\dim(\v)$ and $s=\dim(\z)$ as follows:  Let $\a$ be a one-dimensional inner product
space and $A$ a choice of unit vector in $\a$.  Define an inner product space $(\s, \ip)$
by taking the orthogonal direct sum of $\a$ and $\n$.  Give $\s$ the unique Lie algebra
structure for which
$\n$ is an ideal (the nilradical) and such that
$\ad(A)_{|\v}=\half\Id$ and $\ad(A)_{|\z}= \Id$.  Then $\s$  is a metric solvable Lie
algebra of rank one and eigenvalue type $(1,2;r,s)$.  We will refer to the associated
simply-connected Riemannian solvmanifold $(S,g)$ as the solvmanifold defined by the data
triple $(\v,\z,j)$.

(iv) Data triples $(\v,\z,j)$ and $(\v ',\z ', j')$ will be said to be {\it equivalent}
if there exist orthogonal transformations $\al$ of $\v$ and $\b$ of $\z$ such that 
$$j'(\b(Z))=\al\circ j(Z)\circ\al^{-1}$$ for all $Z\in\z$.
\end{notation}

\begin{remark} The choice of $\half$ and 1 as the eigenvalues of $A$ in Notation
\ref{vz}(iii), as opposed to, say, 1 and 2, is for convenience in the statement of
Theorem \ref{ec} below.  Up to scaling of the metric, all Riemannian solvmanifolds of
algebraic rank one and eigenvalue type $(1,2;r,s)$ arise as in Notation \ref{vz}(iii).  
\end{remark}

\begin{example}\label{heis}
The classical $2n+1$-dimensional Heisenberg algebra is a two-step nilpotent Lie
algebra with basis $\{X_1,\dots, X_n,Y_1,\dots, Y_n, Z\}$ satisfying the
bracket relations $[X_i,Y_i]=Z$ for $i=1,\dots,n$ with all other brackets of basis
elements being zero.  The derived algebra is one-dimensional and is spanned by $Z$. If we
choose the inner product for which the basis above is orthonormal, then $j(Z)X_i=Y_i$ and
$j(Z)Y_i=-X_i$ for all $i$.  The associated solvmanifold is the complex hyperbolic space
of real dimension $2n+2$. 

A two-step nilpotent Lie algebra $\n$, defined as in \ref{vz}(iii) by data
$(\v,\z,j)$,  is said to be of {\it Heisenberg type} if $j(Z)^2=-|Z|^2Id$  for all
$Z\in\z$.  The Heisenberg type algebras, first introduced by A. Kaplan \cite{K}, have
been studied extensively.  Damek and Ricci showed that the solvmanifolds associated as 
in \ref{vz}(iii) to metric nilpotent Lie algebras of Heisenberg type are harmonic (and
thus necessarily Einstein).  We will refer to such solvmanifolds as Damek--Ricci
manifolds.  The fact that these manifolds are Einstein will also follow from Theorem
\ref{ec} below.
\end{example}

\begin{proposition}\label{iso2} The Riemannian solvmanifolds defined, as in
Notation \ref{vz}, by data triples $(\v,\z,j)$ and $(\v ',\z ', j')$ are isometric if and
only if $(\v,\z,j)$ and $(\v ',\z ', j')$ are equivalent.
\end{proposition}

\begin{pf} Let $\s$ and $\s '$ be the metric solvable Lie algebras defined by these data
triples and let $\n$ and $\n '$ be their nilradicals.  By Theorem \ref{iso}, we need only show
that $\s$ and $\s '$ are isomorphic precisely when the data triples are equivalent.  An
elementary argument shows that if the data triples are equivalent, then the metric
Lie algebras $\n$ and $\n '$ are isomorphic.  Moreover, any isomorphism of the metric Lie
algebras $\n$ and $\n '$ extends in an obvious way to an isomorphism of the metric Lie
algebras $\s$ and $\s '$.  Conversely any isomorphism from $\s$ to $\s '$ must restrict
to an isomorphism from $\n$ to $\n '$ which carries $\v$ to $\v '$ and $\z$ to $\z '$. 
Such an isomorphism forces the data triples to be equivalent.  
\end{pf}

\begin{theorem}\label{ec} \cite{EH2} A Riemannian solvmanifold defined by data $(\v,\z,j)$
as in Notation \ref{vz} is Einstein if and only if the following two conditions hold:
\begin{itemize}
\item[(i)]  The map $j:\z\to j(\z)\subset\so(\v)$ is a linear isometry relative to the
Riemannian inner product $\langle\,,\,\rangle$ on $\z$ and the inner product $(\,,\,)$ on
$\so(\v)$ given by $(\alpha,\b)=3D-\frac{1}{r}\tr(\al\b)$.  (Note that $\so(\v)$ is =
isomorphic to $\so(r)$.)
\item[(ii)]  Letting $\{Z_1,\dots ,Z_s\}$ be an orthonormal basis of $\z$,
then $\sum_{i=3D1}^s j(Z_i)^2$ is a scalar operator. (This condition is independent of the
choice of orthonormal basis.)
\end{itemize}
\end{theorem}

\begin{remark} In the presence of condition (i), condition (ii) says that 
${\displaystyle \sum_{i=1}^s j(Z_i)^2 = -s\Id. }$
\end{remark}

\begin{pf} We have $\ad(A)_{|\v}=\half \Id$ and $\ad(A)_{|\z}= \Id$. The mean curvature
vector field $H$, defined in Notation \ref{U}, is given by
$$H=\tfrac 12 (r+2s)A.$$  

Using formula (\ref{ricci}) for the Ricci curvature, the fact that $\a$, $\v$ and $\z$ are orthogonal with respect to the Riemannian inner product, and the fact
that the Killing form $B$ satisfies $B_{|\n \times \s} =0$, one easily checks that the
subspaces $\a$, $\v$ and $\z$ of $\s$ are mutually orthogonal with respect to the Ricci
form (independently of whether the metric is Einstein).  Thus we need only examine the
restriction of the Ricci form to each of these three subspaces.  In applying equation
(\ref{ricci2}), it is convenient to choose an orthonormal basis $\{X_1,\dots,X_n\}$ of
$\s$ so that each $X_i$ lies in one of $\a$, $\v$ or $\z$. 

from equation (\ref{ricci2}) and the facts that $\ad(A)^*=\ad(A)$ and that $A\perp [\s,\s]$, 
we see that 
\begin{equation}\label{a}\ric(A,A)=-B(A,A)=-\tfrac 14 (r+4s).\end{equation}

Next let $X\in\v$.  We have 
$$\ad(X)\circ\ad(X)^*(A)=0=\ad(X)\circ\ad(X)^*_{|\v\cap X^\perp},$$ 
$$\ad(X)\circ\ad(X)^*(X)=\tfrac 14 |X|^2X,$$ 
and
$$\langle\ad(X)\circ\ad(X)^*(Z),Z\rangle =-\langle j(Z)^2X,X\rangle$$ 
for $Z\in\z$ (see equation (\ref{j2})).  Thus
$$-\tfrac 12 \tr\,\ad(X)\circ\ad(X)^*= -\tfrac 18 |X|^2 +\tfrac 12 S(X,X)$$
where 
$$S(X,X)=\langle\sum_{i=1}^s j(Z_i)^2(X),X\rangle.$$
The third term in equation (\ref{ricci2}) equals $\tfrac 18 |X|^2$.  Thus
equation (\ref{ricci2}) yields
\begin{equation}\label{x}\ric(X,X)=-\fourth(r+2s)|X|^2+\half S(X,X).\end{equation}

Finally, let $Z\in\z$. The contribution to the third term in the equation (\ref{ricci2})
for $\ric(Z,Z)$ from those $X_i$'s forming an orthonormal basis of
$\v$ is exactly $\frac{r}{4} |j(Z)|^2$ relative to the inner product
$(\,,\,)$.  The only other non-zero contribution to this term arises from the bracket
of $A$ with $Z/|Z|$.  This latter contribution is exactly cancelled by the first term in
equation (\ref{ricci2}).  Indeed, since
$\image(\ad(Z))=\R Z$ and
$\ad(Z)^*(Z/|Z|)=-|Z| A$, we see that
$\tr(\ad(Z)\circ\ad(Z)^*) =|Z|^2$.  Hence 
\begin{equation}\label{z}\ric(Z,Z)=-\half(r+2s)|Z|^2+\frac{r}{4}|j(Z)|^2. \end{equation}

Assume the metric is Einstein. Equation (\ref{x}) then implies that $S =c_1\langle
\,,\,\rangle _{|\v\times\v}$ for some constant $c_1$ and thus condition (ii) holds.  From
equation (\ref{z}), we see that there exists a constant $c_2$ such that $|j(Z)|^2=
c_2|Z|^2$. Comparing equations (\ref{a}) and (\ref{z}), we see that $c_2=1$.  It then
follows from the definition of $S$ that $c_1=-s$.  Conversely, when conditions (i) and
(ii) hold, we see from equations (\ref{a}), (\ref{x}), and (\ref{z}) that the metric is
Einstein.
\end{pf}

\begin{definition}\label{uniform} We will say an $s$-dimensional subspace $W$ of $\so(r)$
is {\it uniform} if, relative to the inner product on
$\so(r)$ defined in Theorem \ref{ec}, some (and hence every) orthonormal basis
$\{\al_1,\dots \al_s\}$ of $W$ satisfies $\sum_{i=1}^s\alpha_i^2=-c\Id$ for some $c$
(necessarily $c=s$).  We will say two uniform subspaces of $\so(r)$ are {\it
equivalent} if they are conjugate by an orthogonal transformation of $\R^r$.
\end{definition}

\begin{corollary}\label{cor} Up to scaling, isometry classes of Einstein solvmanifolds
of algebraic rank one and eigenvalue type $(1,2;r,s)$ are in one-to-one correspondence
with equivalence classes of uniform $s$-dimensional subspaces of $\so(r)$.  
\end{corollary}

\begin{pf}  If $(S,g)$ is an Einstein solvmanifold with associated data
triple $(\v,\z,j)$, then Theorem \ref{ec} implies that $W:=j(\z)$ is a uniform subspace
of $\so(r)$, where $\so(r)$ is identified with $\so(\v)$.  The equivalence class of $W$
is independent of the choice of isomorphism between $\so(r)$ and $\so(\v)$. Moreover, any
data triple equivalent to
$(\v,\z,j)$ gives rise to a subspace of $\so(r)$ equivalent to $W$.  Conversely, given
a uniform $s$-dimensional subspace of $\so(r)$, define a data triple $(\v,\z,j)$ by
taking $\v$ to be $r$-dimensional Euclidean space, $\z$ to be $W$ equipped with the inner
product $(\,,\,)$ from $\so(\z)$, and $j$ to be the identity map.  Then this data triple
satisfies the two conditions of Theorem \ref{ec} and thus the associated Riemannian
solvmanifold is Einstein.  Equivalent subspaces of $\so(r)$ give rise to equivalent data
triples.  The corollary now follows from Theorem \ref{ec} and Proposition \ref{iso2}.
\end{pf}

\begin{proposition}\label{properties} Let $r\in\Z^+$.
\begin{itemize}
\item[(i)] If $W_1$ and $W_2$ are orthogonal uniform subspaces of $\so(r)$,
then $W_1+W_2$ is also uniform.
\item[(ii)] If $W_1$ and $W_2$ are uniform subspaces of $\so(r)$ with
$W_1\subset W_2$, then $W_2\ominus W_1$ is uniform.
\item[(iii)] If $W$ is a subalgebra of $\so(r)$ and if the restriction to $W$
of the standard representation of $\so(r)$ on $\R^r$ has only one isotypic
component (i.e., all the
irreducible components are equivalent), then $W$ is a uniform subspace of
$\so(r)$.  In particular, $\so(r)$ is uniform in itself.
\item[(iv)] If $W$ is uniform in $\so(r)$, then so is the orthogonal complement
of $W$ in $\so(r)$.
\item[(v)] If $r$ is odd, then $\so(r)$ contains no uniform subspaces of
dimension one or two or of codimension one or two.
\end{itemize}
\end{proposition}

\begin{pf} Statements (i) and (ii) are immediate from Definition \ref{uniform}
and (iv) follows from (ii) and (iii).  To prove (iii), let $W$ be a subalgebra
of $\so(r)$.  The expression $\sum_{i=1}^s\alpha_i^2$ in Definition \ref{uniform} 
is the Casimir operator for the representation of $W$ on $\R^r$.  The Casimir 
operator is constant on each isotypic component of the representation, so (iii) 
follows. 

Finally we prove (v).  Assume $r$ is odd.  It is clear that no one-dimensional
subspace of $\so(r)$ can be uniform in this case and thus, by (iv), no subspace
of codimension one can be uniform.  Next suppose $W$ is a two-dimensional uniform 
subspace of $\so(r)$.  Let $\{a_1,a_2\}$ be an orthonormal basis of $W$.  Since 
$a_1^2+a_2^2$ is a scalar, each $a_i^2$ commutes with $a_1^2+a_2^2$, and thus 
$a_1^2$ commutes with $a_2^2$.  Consequently, $a_2^2$ leaves invariant the 
$0$-eigenspace of $a_1^2$.  Since this eigenspace is odd-dimensional and since
the eigenspaces of $a_2^2$ for non-zero eigenvalues are even-dimensional,
$a_1^2$ and $a_2^2$ must have a common $0$-eigenvector, contradicting
uniformity of $W$.  Thus $\so(r)$ contains no uniform subspaces of dimension 
two or codimension two.
\end{pf}

\begin{example} We classify the Einstein solvmanifolds of algebraic rank one and
eigenvalue type $(1,2;r,s)$ with $r=2$ or $r=3$.  For $r=2$, the algebra $\so(2)$ is 
one-dimensional, so necessarily $s=1$ and $\n$ is the three-dimensional
Heisenberg algebra.  As in Example \ref{heis}, the associated Einstein manifold is
symmetric.  By Proposition  \ref{properties}, for $r=3$, the only uniform subspace
of $\so(r)$ is $\so(r)$ itself.  The associated Einstein solvmanifold is not a
symmetric or Damek--Ricci manifold in the terminology of Example \ref{heis}.  
\end{example}

We next classify the Einstein solvmanifolds of rank one and eigenvalue type
$(1,2;4,s)$; equivalently,  we classify the equivalence classes of uniform subspaces of
arbitrary dimension in $\so(4)$.  The algebra $\so(4)$ is isomorphic to $\so(3)\oplus
\so(3)$.  Viewing  $\R^4$ as the space of quaternions $\bH$, then $\so(4)$ is identified
with $\PP$  where $P$ is the space of purely imaginary quaternions.  The element
$(q,p)\in\PP$  acts on $\H$ as $L(q)+R(p)$ where $L(q)$ and $R(q)$ denote left and right 
multiplication by the quaternion $q$.  The inner product on $\so(4)$ agrees, up to
scalar multiple, with the standard inner product on $\PP$ given by 
$\la(q,p),(q',p')\rangle=q\bar{q'}+p\bar{p'}$.

\begin{lemma}\label{uc}  Suppose $\{\alpha_1,\dots,\alpha_s\}$ is a set of
orthonormal vectors in $\PP=\so(4)$.  Let $B$ be the matrix whose $i^{th}$ 
row is the vector $\alpha_i$ expressed in terms of the standard orthonormal 
basis of $\PP$.  Then the subspace $W$ of $\so(4)$ spanned by $\{\alpha_1,\dots,
\alpha_s\}$ is uniform, if and only if the first three columns of $B$ are orthogonal to 
the last three columns, where the columns are viewed as vectors in $\R^s$ 
with its standard dot product.
\end{lemma}

\begin{pf} For $q, p\in P$, the linear transformation $(q,p)$ of $\bH$ satisfies
\begin{equation}\label{5} (q,p)^2=-(|q|^2 +|p|^2)\Id +2L(q)R(p),
\end{equation}
where $|q|^2=q\bar{q}$.   Letting $\ep_1=\i$, $\ep_2=\j$, and $\ep_3=\k$ be
the standard basis vectors of $P$, write
$\alpha_i=(q,0)+(0,p)=\sum_{j=1}^3[(a_{ij}\ep_j,0)+(0,b_{ij}\ep_j)]$, so
that the $i^{th}$ row of $B$ is given by 
$(a_{i1},a_{i2},a_{i3},b_{i1},b_{i2},b_{i3})$.  Then from equation (\ref{5}) 
we see that
$$\sum_{i=1}^s\alpha_i^2=-\sum_{i=1}^s[|q_i|^2 +|p_i|^2]\Id
+\sum_{j,k=1}^s c_{jk}L(\ep_j)R(\ep_k)$$ where
$$c_{jk}=\sum_{i=1}^s a_{ij}b_{ik}.$$

Note that $c_{jk}$ is the dot product of the $j^{th}$ column of $B$ with the
$(3+k)^{th}$ column.  The lemma follows from the fact that $\sum_{j,k=1}^s
c_{jk}L(\ep_j)R(\ep_k)$ is a multiple of the identity operator only when
all the $c_{jk}$ are zero.
\end{pf}

Recall that the action on $\PP$ by conjugation of the special orthogonal
group $\SO(4)$ consists of pairs of rotations $(R,S)$ sending $(q,p)$ to 
$(R(q),S(p))$.  Conjugation by an orthogonal matrix of determinant $-1$ 
is given by the composition of some such $(R,S)$ with the map that interchanges 
$p$ and $q$.

We will now classify the equivalence classes of uniform subspaces of $\so(4)$
of dimension $s$ for each $s=1,\dots, 6$.  By Proposition \ref{properties}(iv), we
need only classify the uniform subspaces of dimension less than or equal to 3.  Except
when indicated, the examples below are neither symmetric not Damek-Ricci manifolds (see
Example \ref{heis}). 

After conjugating $W$ by an  orthogonal matrix of determinant $-1$ if necessary, we
may assume that the dimension of the subspace of $\R^s$ spanned by the first three
columns of $B$ is greater than or equal to that of the subspace spanned by the last three
columns.

\begin{itemize}
\item[$\bold{s=1}$] In this case we obtain only one solvmanifold, the
6-dimensional complex hyperbolic space.  (Indeed, for any $r$, the only Einstein
solvmanifolds of algebraic rank one and eigenvalue type $(1,2;r,1)$ are the
complex hyperbolic spaces.)

\item[$\bold{s=2}$]  If columns 1--3 of $B$ span a 2-dimensional
space, then columns 4--6 are trivial.  After applying an orthogonal transformation, we
have $\alpha_1=L(\epsilon_1)$ and $\alpha_2=L(\epsilon_2)$.  The resulting metric
nilpotent Lie algebra $\n$ is of Heisenberg type, as defined in Example \ref{heis}, and
the associated Einstein solvmanifold is a Damek-Ricci manifold.  The complementary
four-dimensional uniform subspace is a subalgebra of $\so(4)$ isomorphic to $\su(2)$.

If columns 1--3 span a one-dimensional space, then so do columns 4--6.
Using the fact that exactly two rows of $B$ are independent, after applying an orthogonal
transformation, we can arrange that $\alpha_1=L(\epsilon_1)$ and $\alpha_2=R(\epsilon_1)$.

Up to equivalence, these are the only two possibilities.  Thus we obtain exactly
two isometry classes of Einstein manifolds with $s=2$ and exactly two with $s=4$.

\item[$\bold{s=3}$]  A similar argument shows we get exactly two uniform subspaces up to
equivalence: (i) $W=\{(q,0):q\in P\}$ (the associated solvmanifold is the quaternionic
hyperbolic space);
(ii) $W$ has orthonormal basis $\alpha_1=L(\epsilon_1)$, $\alpha_2=L(\epsilon_2)$ and
$\alpha_3=R(\epsilon_1)$.
\end{itemize}

We have now seen that there are only isolated examples of Einstein manifolds of
eigenvalue type $(\half,1;r,s)$ with $r\leq 4$.  A case-by-case check shows that the same
holds for $r=5$, $s \leq 3$.  In contrast, in the next section we will construct a
continuous family of Einstein manifolds of type $(\half,1;6,3)$.

\section{A Two-Parameter Family}\label{deform}

We will now construct a family of uniform three-dimensional subspaces of $\so(6)$
as follows:  We first construct  three mutually orthogonal three-dimensional uniform
subspaces  $U$, $V$, and $W$ of $\so(6)$ with orthonormal bases $\{A_1,A_2,A_3\}$ of $U$, 
$\{B_1,B_2,B_3\}$ of $V$ and $\{C_1,C_2,C_3\}$ of $W$ such that, for each $j$, 
$A_j$, $B_j$ and $C_j$ anti-commute.  This condition implies that, for any $r,s,t\in\R$, 
$$(rA_i +sB_i+tC_i)^2=r^2A_i^2+s^2B_i^2+t^2C_i^2.$$
Consequently, if $r^2+s^2+t^2=1$, then the set $\{rA_i +sB_i+tC_i: i=1,2,3\}$
is orthonormal and spans a uniform subspace of $\so(6)$.  Of course, $(r,s,t)$
and $(-r,-s,-t)$ determine the same uniform subspace.  We thus obtain a family of uniform
subspaces parameterized by $\R P^2 = S^2/\{\pm1\}$.  After explicitly constructing these
uniform subspaces, we show that some are negatively curved.  Then we test for
equivalences among them.  We can also replace some of the $B_i$ and $C_j$ by their negatives
to get additional families of uniform  subspaces.  

To construct $U$, $V$, and $W$, let $\tau :\fu(3)\to\so(6)$ be the inclusion given
as follows: The algebra $\fu(3)$ consists of all  matrices of the form $X+\sqrt{-1}\,Y$
where $X$ is a skew-symmetric real $3\times 3$ matrix and $Y$ is a symmetric real $3\times
3$ matrix.  Define
$$\tau(X+\sqrt{-1}\,Y)= \pmatrix X & Y \\ -Y & X \endpmatrix.$$  Define the subspaces
$U_0$, $V_0$ and $W_0$ of $\fu(3)$ to be, respectively, the subspace of real
skew-symmetric matrices, the subspace of purely imaginary symmetric matrices with zero
diagonal, and the subspace of purely imaginary diagonal matrices.  Let $U$, $V$, and
$W$ be the images of $U_0$, $V_0$ and $W_0$, respectively, under $\tau$.  Set 
$$X_1=\sqrt{\frac 32}\pmatrix 0 & 0 & 0 \\ 0&0&-1\\0&1&0 \endpmatrix \quad
X_2=\sqrt{\frac 32}\pmatrix 0 & 0 & -1 \\ 0&0&0\\1&0&0 \endpmatrix \text{  and  }
X_3=\sqrt{\frac 32}\pmatrix 0 & -1 & 0 \\ 1&0&0\\0&0&0 \endpmatrix$$
and let $A_i=\tau(X_i)$. Then $\{A_1,A_2,A_3\}$ is
an orthonormal basis of $U$.  Define $Y_1$, $Y_2$ and $Y_3$ analogously to 
$X_1$, $X_2$ and $X_3$ but with all $-1$'s replaced by $1$'s, and let $B_i=
\tau(\sqrt{-1}\,Y_i)$.  Then $\{B_1,B_2,B_3\}$ is an orthonormal  basis of $V$.  For an
orthonormal basis of $W$, let $C_i$ be the image under $\tau$ of the diagonal  matrix
whose $i^{th}$ diagonal entry is $\sqrt{-3}$ and whose other diagonal entries  are zero. 
With these bases in hand, one easily checks that $U$, $V$ and
$W$ are uniform subspaces of $\so(6)$ and that $A_i$, $B_i$ and $C_i$ anti-commute.  

For given $(r,s,t)$ with $r^2+s^2+t^2=1$, let $W(r,s,t) = \text{span}\{D_i = rA_i + sB_i +
tC_i \mid i=1,2,3\}$.  The $D_i$'s give an orthonormal basis for
$W(r,s,t)$.  Since $W(-r,-s,-t)=W(r,s,t)$, we have a two-dimensional moduli space
of uniform subspaces of $\so(6)$,  parametrized by $\R P^2$.  Each uniform subspace gives
rise to an Einstein  solvmanifold as in Corollary \ref{cor}.  
%
%
%
%
%

\begin{proposition}\label{sect1} The deformation includes a continuous family of
Einstein solvmanifolds with negative sectional curvature.
\end{proposition}

\begin{pf} We show that the Einstein solvmanifold corresponding to the uniform
subspace $W(1,0,0)$ has negative sectional curvature.  By continuity, $W(r,s,t)$ must also
have negative curvature for all $(r,s,t)$ sufficiently close to $(1,0,0)$. 

Recall that the sectional curvature of the two-plane in $\s$ spanned by
orthonormal vectors $X$ and $Y$ is given by 
$$K(X,Y) = -\frac 34 |[X,Y]|^2 - \frac 12 \langle [X,[X,Y]],Y\rangle -
\frac 12 \langle [Y,[Y,X]],X\rangle + |U(X,Y)|^2 - \langle U(X,X),U(Y,Y) \rangle$$
where $U$ is defined as in Notation \ref{U}.  (See Besse \cite{Bes}, equation (7.30).)

By Proposition 2.3.1 of Leukert \cite{L}, in order to show that $S$ has negative sectional
curvature, it suffices to show that $K(P,Q) < -\frac 34 |[P,Q]|^2$, whenever $P$ and $Q$ are
orthonormal vectors in the nilradical $\n$. By rotating the basis $\{P,Q\}$ within the plane
spanned by these vectors, we obtain a new orthonormal basis $\{X+Z,Y+W\}$, where $X$ and $Y$
are orthogonal vectors in $\v=\R^6$ and $Z$ and $W$ are orthogonal vectors in $\z$.  We
can then  rewrite the sectional curvature as 
$$K(X+Z,Y+W)=-\frac 34|[X,Y]|^2 -(\frac 12 |X|^2+|Z|^2)(\frac
12|Y|^2+|W|^2) -\langle j(Z)X,j(W)Y\rangle +\frac 14|j(Z)Y+j(W)X|^2.$$ 

To see that $K(X+Z,Y+W) < -\frac 34 |[X,Y]|^2$, we need only show that 
\begin{equation}\label{show}(\frac 12 |X|^2+|Z|^2)(\frac 12|Y|^2+|W|^2) 
> -\langle j(Z)X,j(W)Y\rangle  +\frac 14|j(Z)Y+j(W)X|^2.\end{equation} 

The uniform subspace $W(1,0,0)$ of $\so(6)$ was obtained by diagonally embedding $\so(3)$ 
into $\so(6)$. In \cite{L}, Leukert proved that the seven-dimensional Riemannian
solvmanifold constructed from $\so(3)$ (as a uniform subspace of itself) is negatively curved
by verifying inequality (\ref{show}) in this case.  (The metric on his solvmanifold was
rescaled but this does not affect the computations.)  It is straightforward to verify that
inequality (\ref{show}) in our case follows from the case proven by Leukert.
\end{pf}
%
%
%

We next verify that this deformation of Einstein manifolds is non-trivial:

\begin{proposition}\label{conj} The set of isometry classes of the Einstein
manifolds $W(r,s,t)$ defined above is parameterized by the quotient of $\R P^2$ by a
finite equivalence relation. 
\end{proposition}

\begin{lemma}\label{centralizer} Let $W=W(r,s,t)$ be a uniform subspace of $\so(6)$
as above.  The centralizer of $W$ is one-dimensional, spanned by $C_1+C_2+C_3$.
\end{lemma}

\begin{pf} The eigenvalues of $D_i^2$ are $-t^2$ and $-(r^2+s^2)/2$. First assume $t^2
\neq (r^2+s^2)/2$. The eigenspace of $D_i^2$ corresponding to $-t^2$ is the non-zero
eigenspace of $C_i$.   For $i=1$, this is $\{\bold{e_1, e_4}\}$; for $i=2$, it is
$\{\bold{e_2, e_5}\}$; for $i=3$,  it is $\{\bold{e_3, e_6}\}$, where $\{\bold{e_1, \dots
e_6}\}$ is the standard basis of $\R^6$.   Any element $X$ in $\so(6)$  which commutes
with all of $W$ must leave each of three subspaces invariant.  It then follows that $X$
is a linear combination of
$C_1,\, C_2,\, C_3$.   It is easy to check that the only combinations which commute with
all of $W$ are those for  which all coefficients are equal.

When $t^2=(r^2+s^2)/2$, the centralizer of $W$ is still spanned by $C_1+C_2+C_3$.  
This is easily checked (e.g. using Maple).
\end{pf}

\begin{lemma}\label{angle1} Let $W(r,s,t)$ be a uniform subspace of $\so(6)$ as above.  
The angle $\theta$ between $W(r,s,t)$ and its centralizer satisfies $\cos\theta = |t|$.
\end{lemma}

\begin{pf} By Lemma \ref{centralizer},  $\theta$ is the angle 
between $W$ and $C_1+C_2+C_3$.  Any unit vector in $W(r,s,t)$ is of the form 
$xD_1+yD_2+zD_3$, for $x,y,z \in \R$ satisfying $x^2+y^2+z^2=1$. Since $|C_1+C_2+C_3|=
\sqrt 3$, we have
$$\cos\theta = \max \frac{\langle C_1+C_2+C_3, xD_1+yD_2+zD_3 \rangle}
{|C_1+C_2+C_3|\cdot |xD_1+yD_2+zD_3|} =\max\frac{|t(x+y+z)|}{\sqrt 3}=|t|.$$
\end{pf}

\begin{lemma}\label{relns} Let $W(r,s,t)$ be a uniform subspace of $\so(6)$ as above.  
Then the minimal angle $\theta$ between $W(r,s,t)$ and $[W(r,s,t),W(r,s,t)]$ satisfies
$$\cos\theta = \frac{|r|\sqrt{1-t^2}}{\sqrt{1+3t^2-2|t^2+\sqrt 2\,st|}}.$$
\end{lemma}

\begin{pf}
We compute the generating elements of $[W(r,s,t),W(r,s,t)]$:
\begin{align}
[D_1,D_2] &= \frac{r^2+s^2}{\sqrt 2}A_3 + t(r(B_2-B_1)+s(A_1-A_2)) \\
[D_3,D_1] &= \frac{r^2-s^2}{\sqrt 2}A_2 +\sqrt{2}rsB_2 +
t(-r(B_1+B_3)+s(A_1+A_3)) \\
[D_2,D_3] &= \frac{r^2+s^2}{\sqrt 2}A_1 + t(r(B_2-B_3)+s(A_3-A_2)) 
\end{align}
We find the non-zero inner products of these basis elements of
$[W(r,s,t),W(r,s,t)]$ with basis elements of $W(r,s,t)$:
$$\langle [D_i,D_j],D_k \rangle = \frac{r(r^2+s^2)}{\sqrt 2} \quad
\text{for $(i,j,k)$ cyclic permutations of $(1,2,3)$}.$$ 

Let $x,y,z,u,v,w \in \R$ satisfy $x^2+y^2+z^2=1$ and $u^2+v^2+w^2=1$.  Then
\begin{equation}
\langle x[D_1,D_2] +y[D_3,D_1] +z[D_2,D_3], uD_1 +vD_2 +wD_3 \rangle 
= (xw+yv+zu)\frac{r(r^2+s^2)}{\sqrt 2}.
\end{equation}
Since $|uD_1 +vD_2 +wD_3 | = 1$ and  
$$|x[D_1,D_2] +y[D_3,D_1] +z[D_2,D_3] | = 
\sqrt{(r^2+s^2)( \tfrac{r^2+s^2+4t^2}{2} +2(xz+xy+zy)(t^2 +\sqrt 2 st))},$$
we have 
\begin{equation}
\cos\theta =\max \frac{|(xw+yv+zu)r|\sqrt{r^2+s^2}}{\sqrt{(r^2+s^2+4t^2) +
4(xz+xy+yz)(t^2+\sqrt 2 st)}}.
\end{equation}
The maximum is attained when $(w,v,u)=(x,y,z)$ and
$xz+xy+yz=\pm 1/2$ with the sign chosen so that $(xz+xy+yz)(t^2+\sqrt 2 st) \leq 0$.  
This gives 
\begin{equation}\label{brang}
\cos\theta =|r| \sqrt{\frac{r^2+s^2}{r^2+s^2+4t^2-2|t^2+\sqrt 2 st|}}.
\end{equation} 
Replacing $r^2+s^2$ by $1-t^2$, we obtain the expression in the lemma.
\end{pf}

We can now prove Proposition \ref{conj}.
 Suppose $W(r',s',t')$ is equivalent to $W(r,s,t)$.  From 
Lemma \ref{centralizer}, we know that $W(r,s,t)$ and $W(r',s',t')$ share the
same centralizer.  Since the angle  between the uniform space and its centralizer is
preserved under conjugation by elements of $\O(6)$, it follows from Lemma \ref{angle1}
that $|t| = |t'|$ and thus also
$r^2+s^2=(r')^{2}+(s')^{2}$.   Since the angle in Lemma \ref{relns} is also
preserved under conjugation, it follows that $r'$ and $s'$ are determined up to finitely
many possibilities.

We remark that the finite equivalence relation in Proposition \ref{conj} is non-trivial.
For example the element of $\O(6)$ with the $3\times 3$ identity matrix in the upper
right and lower left $3\times 3$ blocks and zeros elsewhere conjugates $W(r,s,t)$ to
$W(r,-s,-t)$.
%
%
%
%
%
\section{Examples: Irreducible Symmetric Spaces of Type IV}
We now obtain new examples of Einstein manifolds by modifying the Lie algebra
structure of the solvable Lie algebras associated with the non-compact
symmetric spaces of rank greater than one.  We first recall the structure of
the isometry groups of the symmetric spaces of non-compact type.  See Helgason
\cite{H} for more details.

\begin{notation}\label{note1} The identity component of the full isometry
group of an irreducible symmetric space $M$  of non-compact type is a simple Lie
group
$G_0$, and the isotropy subgroup $K$ at a point of $M$ is a maximal compact subgroup of
$G_0$.  The Lie algebra $\rg$ of $G_0$ admits an Iwasawa decomposition $\rg=\fk +\a
+\n$ where $\fk$ is the Lie algebra of $K$, $\a$ is abelian, $\n$ is nilpotent and
$\s:=\a+\n$ is solvable with nilradical $\n$.  The Lie subgroup $S$ of $G_0$ with Lie
algebra $\s$ acts simply transitively on $M$, and $M$ may be identified with
$S$, endowed with a left-invariant metric.

The elements of $\ad_\s(\a)$ are symmetric and $\n$ decomposes into root spaces
$\n_\alpha$, $\alpha\in\S$, where $\S$, the set of positive roots of $\a$
in $\g_0$, is a subset of the dual space $\a^*$ of $\a$.   Relative to the inner product
on $\s$ defined by the symmetric space metric on $S$ (identified with $M$), the root
spaces are mutually orthogonal and there exists an orthonormal basis of the subspace $\n$
of $\s$ such that:
\begin{itemize}
\item Each basis vector lies in some $\n_\b$.
\item The bracket of any two basis vectors is a scalar multiple of another
basis vector.
\item If $X$, $Y$, and $U$ are vectors in the basis with $Y\neq U$, then
$[X,Y]\perp [X,U].$
\end{itemize}
\end{notation}

\begin{notation}\label{note2} Let $\g$ be the complexification of $\g_0$ and let $\g^\R$
be $\g$ viewed as a real Lie algebra.  The dimension of $\g^\R$ is twice that of $\g$. 
View $\g_0$, and thus also $\s$, as  subalgebras of $\g^\R$.  Choose a basis ${\cal B}$ for
$\n$ satisfying the three conditions above.  Consider a new subspace of $\gr$ with basis
${\cal B}'$, where ${\cal B}'$ is obtained from ${\cal B}$ by replacing some of the
vectors $X\in {\cal B}$ by $\sqrt{-1}\,X$.  If the new subspace $\n '$ is closed under the
bracket operation in $\g$, then $\n '$ will be a nilpotent algebra. The brackets
of the basis vectors in $\n '$ just differ by a sign from the corresponding brackets in
$\n$.  Endow $\n '$ with the inner product for which ${\cal B}'$ is orthonormal.  The
subalgebra $\a$ of $\g^\R$ normalizes $\n '$; indeed, $\n '$ decomposes into root spaces
$\n '_\alpha$ of $\a$, where $\alpha$ varies over $\S$.   Letting $\s '=\a+\n '$, then 
$\s'$ is a solvable Lie algebra with nilradical $\n '$.  We will say $\s$ and
$\s '$ are {\it associated} subalgebras of $\g^\R$. Give $\s '$ the inner product for
which $\s '=\a+\n '$ is the orthogonal direct sum of the inner product spaces $\a$ and 
$\n'$.  This inner product on $\s '$ defines a left-invariant Riemannian metric on the
associated simply connected Lie group $S'$. The Riemannian manifold $S'$ will be referred
to as a {\it Riemannian solvmanifold associated with the symmetric space} $M$.
\end{notation}

\begin{theorem}\label{thm1} If $S'$ is a Riemannian solvmanifold associated with
a symmetric space $M$ as in Notation \ref{note2}, then $S'$ is Einstein with the same
Einstein constant as $M$.
\end{theorem}

\begin{pf} We write $\ric$ for the Ricci tensor on the symmetric solvmanifold
$S$, identified with $M$,  and $\ric '$ for the Ricci tensor of $S'$.  View $\ric$ and
$\ric '$ as quadratic forms on $\s$ and $\s '$, respectively.  We denote by $\tau$
the linear isomorphism from $\s$ to $\s '$ which restricts to the identity on $\a$ and 
sends the basis ${\cal B}$ of $\n$ to the basis ${\cal B}'$ of $\n '$ in the obvious
way. Using the expression  for the Ricci curvature given in equation (\ref{ricci}), we
will show that $\ric (X,Y)=\ric '(\tau(X),\tau(Y))$ for all $X$, $Y\in\s$.
 
Let $H$ be the mean curvature vector of $\s$ as in Notation \ref{U}. From the fact that
$[A,\tau(X)]=\tau([A,X])$ for all $A\in\a$ and $X\in\s$, we see that 
$\tau(H)=H$ is also the mean curvature vector of $\s'$ and that $\langle U(X,Y),H
\rangle =\langle U(\tau(X),\tau(Y)),H\rangle$ for all $X$, $Y\in\s$, where $U$ is defined
in Notation \ref{U}.  Also $\tau$ is an isometry between  the Killing forms
$B$ of $\s$ and $B'$ of $\s'$ since the two Killing forms agree on $\a$ and are zero on
the nilradicals $\n$ and $\n'$, respectively. It follows that $\ric(X,Y)=\ric'
(\tau(X),\tau(Y))$ whenever $X\in\a$ and $Y\in \s$. Moreover, from the properties of
${\cal B}$ indicated by bullets, we see that $\ric (X,Y)=\ric' (\tau(X),\tau(Y))$ when
$X$, $Y\in {\cal B}$. \end{pf}

Before applying this theorem, we discuss conditions under which two
Riemannian solvmanifolds associated to a symmetric space are isometric.  

\begin{proposition}\label{isoroot}  In the notation of Notation \ref{note1}
and \ref{note2}, suppose that $S'$ and $S''$ are two Riemannian solvmanifolds associated
with $M$.  Then $S'$ is isometric to $S''$ if and only if there exists an automorphism
$\phi:\S\to\S$ of the root system $\S$ (i.e., a Weyl group element) and an isomorphism
$\tau$ between the metric Lie algebras $\n'$ and $\n''$ such that
$\tau(\n'_\al)=\n''_{\phi(\al)}$ for all $\al\in\S$.  
\end{proposition}

\begin{pf} By Proposition \ref{iso}, $S'$ and $S''$ are isometric if and only if
there exists a metric Lie algebra isomorphism $\tau:\s '\to\s ''$.  Such an isomorphism
$\tau$  must carry $\n '$ to $\n ''$ and thus $\a$ to $\a $.  The restriction of $\tau$
to $\a$ must then induce an automorphism of $\S$ and the restriction of $\tau$ to $\n'$
must satisfy the last statement of the proposition.  
\end{pf}

\begin{example}\label{ex1}  Suppose $\s '$ and $\s ''$ are two solvable subalgebras of
$\g^\R$ associated with $\s$.  If for each $\alpha\in\S$, we have either $\n ''_\al 
=\n'_\al $ or $\n ''_\al =\sqrt{-1}\,\n '_\al $, then the associated Riemannian
solvmanifolds are isometric.  To see this, choose a base $\{\al_1, \dots , \al_k\}$ for
the root system $\Sigma$. Write $\n_j$ for $\n_{\al_j}$. Every root $\al\in\S$
can be written uniquely as a linear combination of  $\{\al_1, \dots ,
\al_k\}$  with non-negative integer coefficients.  Since $[\n_\al ,\n_\b ]
=\n_{\al +\b}$ whenever $\al$, $\b$, $\al +\b\in\S$, the $\n_j$ generate
$\n$ as a Lie algebra. A similar statement holds in the algebras $\n '$ and $\n ''$.

Let
$${\cal A}=\{j\in\{1,\dots,k\}: \n ''_j =\sqrt{-1}\,\n '_j \}.$$
For $\alpha\in\Sigma$, expressed as a linear combination $\al
=\Sigma_{i=1}^k b_i\al_i$, define
the {\it restricted height} $rh(\al)$ by
$$rh(\al)=\Sigma_{j\in{\cal A}}b_j.$$
The fact that $\n '$ and $\n ''$ are both closed under the bracket
operation in $\g^\R$
then forces $\n ''_\al$ to satisfy
$$\n_\al ''=\sqrt{-1}^{rh(\al)}\n '_\al.$$  The linear map from $\n '$ to
$\n ''$ whose
restriction to each $\n '_\al$ is given by multiplication by
$\sqrt{-1}^{rh(\al)}$ is then a Lie
algebra isomorphism and inner product space isometry.
This map extends to an isomorphism from $\s '$ to $\s ''$ by acting as the
identity on $\a$.
By Proposition
\ref{iso}, we conclude that $S''$ is isometric to $S'$.  In the notation of
Proposition \ref{isoroot}, $\phi$ is the identity.
\end{example}

\begin{remark}\label{rhremark} Example \ref{ex1} shows that given any subalgebra $\s '$
of $\g^\R$ associated with $\s$ and any subset ${\cal A}$ of $\{1,\dots,k\}$, where $k$
is the rank of $\s$, then we can define another subalgebra $\s ''$ of $\g^\R$ associated
to $\s $ by setting $\n_\al ''=\sqrt{-1}^{rh(\al)}\n '_\al$ in the notation of
the example.  The associated Riemannian solvmanifolds $S'$ and $S''$ are then isometric. 
This construction will be useful in uniqueness arguments below.
\end{remark}

We now apply Theorem \ref{thm1}.  The irreducible symmetric spaces of
negative Ricci curvature may be divided into two types by considering the structure of
the complex Lie algebra $\g$.  We have either:
\begin{itemize}
\item[(i)] $\g$ is a complex simple Lie algebra, or
\item[(ii)] $\g$ is isomorphic to the direct sum of two copies of a simple
complex Lie algebra.
\end{itemize}
Helgason \cite{H} calls these two types of symmetric spaces ``Types III and IV'',
respectively, with types I and II referring to symmetric spaces of positive Ricci
curvature.  The book \cite{H} is an excellent source for details on the facts about the
root systems of semisimple Lie algebras used below.

We first consider case (ii); i.e., we assume that $\g$ consists of two
copies of a
complex simple Lie algebra $\f$.  In this case, $\g_0$ is isomorphic to
$\f^\R$. Let $J$ denote
the complex structure on $\g_0$ arising from $\f$. When a subspace $V$ of
$\g_0$ is invariant
under $J$, we will denote by $V_{\bC}$ the complex subspace of $\f$ defined
by $V$.  In the
notation above, the subspace $\a$ of $\g_0$ is invariant under $J$, and
$\a_{\bC}$ is a Cartan
subalgebra of $\f$.  Moreover, each root $\al\in\S$ commutes with $J$ and
thus defines a complex
linear functional, again denoted $\al$, on $\a_{\bC}$.  With this
identification and with a
suitable choice of ordering on the set of roots of $\a_\bC$ on $\f$, the set
$\S$ forms the system of positive roots of $\a_{\bC}$ in $\f$, and for
each $\al$, we have
$\f_\al=(\n_\al)_{\bC}$.  In particular, since the root spaces $\f_\al$
have complex dimension
one, the root spaces
$\n_\al$ have real dimension two.  The almost complex structure $J$ on $\s$
is skew-symmetric
relative to the inner product defined by the symmetric space metric on $S$.

One may choose a unit vector $X_\al$ in each $\n_\al$, $\al\in\S$, so that
$$[X_\al,X_\b]=N_{\al,\b}X_{\al+\b}$$
for some non-zero constants $N_{\al,\b}$ whenever $\al,\b\in\S$, with the
convention that
$X_{\al+\b}=0$ when $\al+\b$ is not a root.  The collection of vectors
$${\cal B}_0=\{X_\al, JX_\al:\al\in\S\}$$
forms an orthonormal basis for $\n$.  The brackets satisfy
\begin{equation}\label{brackets}[X_\al,X_\b]=N_{\al,\b}X_{\al+\b}=-[JX_\al,JX_\b
]
\quad \text{and} \quad
[X_\al,JX_\b]=N_{\al,\b}JX_{\al+\b}=[JX_\al,X_\b].\end{equation}

\begin{theorem}\label{thm2}  Let $M$ be a symmetric space of Type IV; i.e., assume its isometry
algebra satisfies
condition (ii).  If $\g$ has rank greater than one, then, up to
isometry, there exists a
unique non-symmetric Einstein solvmanifold $S'$ associated with $M$.  Its
Lie algebra is
spanned by $\a$ together with $\{X_\al,\,JX_\al:\al\in\S\}$.  If
$\g$ has rank one,
there are no non-symmetric Einstein manifolds associated with $M$.
\end{theorem}

\begin{pf}  If $\g$ has rank one, then $\n$ is abelian.  Any solvable
algebra associated
with
$\s '$ will also have abelian nilradical, and it follows from Proposition
\ref{isoroot} that the
associated solvmanifold will be isometric to the symmetric space.

Now assume that $\g$ has rank greater than one. The fact that $\s '$, as
defined in the theorem,
is a solvable subalgebra of $\g$ associated with $\s$ is immediate from the
bracket relations
(\ref{brackets}).  Thus Theorem \ref{thm1} guarantees that $S'$ is Einstein.  The fact
that $S'$ is not
isometric to $S$ can be seen either from Proposition \ref{iso} or by
noting that $S'$ has some
positive sectional curvature.
To see the latter fact, recall that the sectional curvature of the two plane in $\s '$ spanned by a pair of orthonormal vectors $X$ and $Y$ is given by
\begin{equation}\label{sect}K(X,Y) = -\frac 34 |[X,Y]|^2 - \frac 12 \langle
[X,[X,Y]],Y\rangle -
\frac 12
\langle [Y,[Y,X]],X\rangle + |U(X,Y)|^2 - \langle U(X,X),U(Y,Y) \rangle
\end{equation}
where $U$ is defined as in Notation \ref{U}.  (See Besse \cite{Bes}, equation 7.30.) 
Since $\n$ is not abelian, we can choose two basic roots $\al$ and $\b$ in $\S$ such that
$\al+\b$ is a root. Let $X= X_\al+\,JX_\al$ and $Y=X_\b-\,JX_\b$. By
equation (\ref{brackets}), $[X,Y]=0$.  Moreover,
$U(X,Y) = 0$ since $\al$ and $\b$ are basic roots.
Thus the only non-zero term in the right-hand-side of equation (\ref{sect}) is the last
term, and we need to show that 
\begin{equation}\label{uxy}\langle U(X,X),U(Y,Y) \rangle <0.
\end{equation} 

Equation (\ref{uxy}) follows from the structure theory of semisimple Lie algebras
of non-compact type.  Given a root vector $\delta\in\S$, let
$H_\delta\in\a$ be the unique vector satisfying
$\delta(Z)=\langle Z,H_\delta\rangle$ for all
$Z\in \a$.  Then $\langle H_\delta, H_\epsilon\rangle \leq 0$ with strict inequality
holding when $\delta +\epsilon\in\S$. (See \cite{H}.)  In our situation, we have 
$U(X,X)=|X|^2 H_\al$ since $[Z,X]=\al(Z)X$ for all $Z\in\a$, and similarly
$U(Y,Y)= |Y|^2 H_\b$.  
Thus equation (\ref{uxy}) follows and we conclude that $K(X,Y)>0$.

We next turn to the uniqueness statement.

The construction in Notation \ref{note2} of solvmanifolds associated with $M$
depends a priori on a
choice of orthonormal basis
$\cal B$ of $\n$ satisfying the conditions in Notation \ref{note1}.  We
first show under the
hypotheses of Theorem \ref{thm2} that all solvmanifolds associated as in Notation \ref{note2} with $M$ are isometric to ones constructed using the particular
basis ${\cal B}_0$ given above.  Let $\cal C$ be
another orthonormal
basis of $\n$ satisfying the conditions in Notation \ref{note1}.  $\cal C$
contains an
orthonormal basis of each root space $\n_\al$.  This basis may be expressed
in the form $\{U_\al,
JU_\al\}$ for some vectors
$U_\al\in\n_\al$.  Within the complex algebra $\g$, we can write
$U_\al=e^{i\theta_\al}X_{\al}$
for some
$\theta_\al\in[0,2\pi)$.  We then have
$$[U_\al,U_\b]= e^{i(\theta_\al+\theta_\b
-\theta_{\al+\b})}N_{\al,\b}U_{\al+\b}.$$
The constants of structure of $\n$ are real, so
$e^{i(\theta_\al+\theta_\b-\theta_{\al+\b})}
=\pm 1$.  Choose a base $\{\al_1, \dots , \al_k\}$ for the
root system $\Sigma$.  Each root $\al$ can be written uniquely in the form
$\Sigma_{1\leq j\leq k}b_j\al_j$.  Define a linear isomorphism $\tau$ from
$\n$ to itself
by sending $X_\al$ to ${\epsilon_\al}U_\al$ and $JX_\al$ to
${\epsilon_\al}JU_\al$ where $\epsilon_\al=\pm 1$ is given by
$\epsilon_\al=e^{i((\Sigma_{1\leq
j\leq k}b_j\theta_{\al_j})-\theta_{\al})}$.  Then $\tau$ is easily seen to be a metric Lie algebra
isomorphism.  In particular, the basis ${\cal C}'=\{\epsilon_\al U_\al,\,\epsilon_\al JU_\al:\,\al\in\S\}$
satisfies exactly the same
bracket relations as ${\cal B}_0$ and thus any solvable metric Lie algebra $\s ''$ associated with $\s$ which is constructed using ${\cal C}'$ is isomorphic to one
constructed using
$\cal B$ and conversely.    Finally, since
${\cal C}'$ can be obtained from
${\cal C}$ simply by multiplying some elements of
$\cal C$ by
$-1$, the metric Lie algebras associated to $\s$ which can be constructed using
$\cal C$ are precisely the ones that can be constructed using ${\cal C}'$.  Consequently, up to isometry, our construction of Riemannian solvmanifolds associated to $M$ is independent
of the choice of basis satisfying the conditions of Notation \ref{note1}, and we can restrict attention to the basis ${\cal B}_0$.

Thus suppose $\s ''$ is a solvable subalgebra of $\g^\R$ constructed as in
Notation \ref{note2}
using the basis
${\cal B}_0$ of $\n $. We need to show that the associated Riemannian
solvmanifold
is isometric either to the symmetric space $M$ or to the solvmanifold $S'$
of Theorem
\ref{thm2}.  Choose a base $\{\al_1, \dots , \al_k\}$ for the
root system $\Sigma$ and write $\n_j$ for $\n_{\al_j}$.  Let
$${\cal A}=\{j\in\{1,\dots,k\}: \n_j ''\neq\n_j\}.$$
In view of Remark \ref{rhremark}, for $j\in{\cal A}$ we may assume
that $\n_j ''$ is
obtained from
$\n_j$ by multiplying one of
$X_j$ or $JX_j$, but not both, by $\sqrt{-1}$, where $X_j$ denotes
$X_{\al_j}$.  Moreover, by
further multiplying the eigenspace by $\sqrt{-1}$ if necessary and again
referring to Remark \ref{rhremark}, we may assume that
$X_j ''=X_j$ for all
$j$, that $JX_j ''=\,X_j$ when $j\in{\cal A}$, and that $JX_j ''=JX_j$
otherwise.  But if
$j\in{\cal A}$, $i \not\in {\cal A}$, and $\al_j+\al_i\in\S$, then from
the bracket relations (\ref{brackets}), we see that it is impossible to define
$\n_{\al_j+\al_i}''$ in such a way that $\n ''$ is closed under the bracket
operation.  We
conclude that $\S$ decomposes into the direct sum of two root systems, one
with base
$\{\al_j:j\in\cal A\}$ and the other with base
$\{\al_j:j \not\in {\cal A}\}$.  This contradicts the irreducibility of the
symmetric space $M$,
unless $\cal A$ is either empty or all of $\{1,\dots,k\}$.  These two cases
give, respectively,
the symmetric space or the manifold $S'$ defined in Theorem \ref{thm2}.
\end{pf}

 We next consider the case (i); i.e., we assume that  $\g$ is simple.  (See Notation \ref{note2}.)
Letting $\m_0$ be the
centralizer of $\a$ in $\fk$, then the complexification $\h$ of $\a +\m_0$
is a Cartan subalgebra
of $\g$.  The Lie algebra $\g$ decomposes into root spaces
$$\g=\Sigma_{\b\in\Delta}\g_\b$$
where $\Delta$  is a subset of $\h^*$.  For $\b\in\Delta$, let $\bb$ denote
the restriction of
$\b$ to $\a$.  Then
$$\n=\g_0\cap\Sigma_{\bar{\b}\in\S}\g_\b.$$
Each $\g_\b$ is a one-dimensional complex subspace of $\g$. 
However, the $\n_\al$
may have dimension greater than one, since there may be several roots $\b$
such that $\bb=\al$.
The elements $\ad(X)$, $X\in\m_0$ normalize each root space $\n_\al$.  In
case $\m_0$ is
trivial, i.e., $\h$ is the complexification of $\a$, then $\g_0$ is said to
be a {\it normal real
form} of $\g$.  In this case, all the root spaces $\n_\al$ are
one-dimensional.  If $\m_0$ is
non-trivial, then some of the root spaces will have higher dimension.  

\begin{proposition}\label{prop2}If each $\n_\al$ is one-dimensional, i.e.,
if $\g_0$ is a
normal real form of $\g$, then any Einstein manifold associated to $M$ is
isometric to $M$.
\end{proposition}

The proposition is an immediate consequence of Example \ref{ex1}.

We now show by specific construction that in most cases of classical irreducible
symmetric spaces $M$ for which $\g_0$ is not a normal real form of $\g$, there exist
non-symmetric Einstein manifolds associated with $M$.

%
%
%
%
\setcounter{subsection}{8}
Recall that the classical symmetric spaces are given as follows:
\begin{center}
\begin{tabular}{|c|c|c|}
Symmetric Space $M = G/H$ & Rank($M$) & Rank($G$) \\
\hline
$\SO(p,q)/\SO(p)\SO(q)$ & min$(p,q)$ & $\lfloor \frac{p+q}{2}\rfloor$\\
$\SU(p,q)/S(\U(p)\U(q))$ & min$(p,q)$ & $p+q-1$\\
$\Sp(p,q)/\Sp(p)\Sp(q)$ & min$(p,q)$ & $p+q$\\
$\SO(n,\bH)/\U(n)$ & $\lfloor \frac n2 \rfloor$ & $n$\\
$\SL(n,\bH)/\Sp(n)$ & $n-1$ & $2n-1$\\
$\SL(n,\R)/\SO(n)$ & $n-1$ & $n-1$\\
$\Sp(n,\R)/\U(n)$ &$n$ & $n$ \\
\hline
\end{tabular}
\end{center}

We note that Helgason\cite{H} refers to $\SO(n,\bH)$ and $\SL(n,\bH)$ as
$\SO^*(2n)$ and $\SU^*(2n)$, respectively.

The Lie algebra $\so(p,q)$ is a real form of $\so(n,\bC)$ where $n=p+q$.  Since
$\SO(p,q)$ is isomorphic to $\SO(q,p)$, we may assume $p\leq q$.  In case $q=p$
or $q=p+1$, then $\so(p,q)$ is the normal real form of $\so(n,\bC)$.  In all other
cases, i.e., whenever $q-p\geq 2$, we will construct in Subsection \ref{SO(p,q)}
non-symmetric Einstein solvmanifolds associated with the symmetric space
$M=\SO(p,q)/\SO(p)\SO(q)$ as in Notation \ref{note2}.

The Lie algebras $\su(p,q)$ are real forms of $\sl(p+q,\bC)$.  None of these are
normal forms; the normal real form of $\sl(n,\bC)$ is $\sl(n,\R)$.  However, we
will obtain in Subsection \ref{SO(p,q)} non-symmetric Einstein solvmanifolds
associated with the symmetric space $M=\SU(p,q)/S(\U(p)\U(q))$ only when $q\geq
p+2$.  (As in the case of $\so(p,q)$, we may take $q\geq p$.)  Similarly,
$\fsp(p,q)$ is a real form of $\fsp(n,\bC)$, while $\fsp(n,\R)$ is the normal real
form, but we will obtain new Einstein solvmanifolds only when $q\geq p+2$.  

Finally, we will obtain non-symmetric Einstein manifolds associated with
$\SO(n,\bH)/\U(n)$ in Subsection \ref{SO(n,H)} and with $\SL(n,\bH)/\Sp(n)$ in
Subsection \ref{SL(n,H)}.

\subsection{SO(p,q)/SO(p)SO(q), SU(p,q)/S(U(p)U(q)), and 
Sp(p,q)/Sp(p)Sp(q)}\label{SO(p,q)}
We consider simultaneously the rank $p$ non-compact
irreducible symmetric spaces
$\SO(p,q)/S(\O(p)\O(q))$,\\ $\SU(p,q)/S(\U(p)\U(q))$ and
$\Sp(p,q)/\Sp(p)\Sp(q)$ dual to the real, complex and quaternionic
Grassmannians.  We assume that $q-p\geq 2$.  The Lie algebras $\so(p,q)$,
$\su(p,q)$ and $\fsp(p,q)$ are given by
$$\g_0 = \{ A \in \gl(p+q,\bF) \mid AM + M\bar A^t = 0 \} \text{ where }
M = \pmatrix \Id_p & 0 \\ 0 & -\Id_q \endpmatrix$$
where   $\bF=\R$ in the case of $\so(p,q)$, $\bF=\bC$ in the case of
$\su(p,q)$ and $\bF=\bH$ in the case of $\fsp(p,q)$.

Let $m=q-p$.  Then
$$\g_0=\pmatrix A & B & C \\  \bar{B}^t & D & E \\ \bar{C}^t & -\bar{E}^t
& F \endpmatrix$$
where $A$, $B$, and $D$ are $p\times p$ matrices, $C$ and $E$ are $p\times
m$ matrices and $F$ is an $m\times m$ matrix.  The matrix $A$ satisfies
$\bar{A}^t+ A =0$ and similarly for $D$ and $F$; the matrices $B$, $C$ and $E$ are
arbitrary.

The Cartan subalgebra $\a$ consists of those matrices for which $B$ is
real and diagonal and all other blocks are zero.  Letting
$\{\om_1,\dots,\om_p\}$ be the basis of linear functionals on $\a$ dual to
the standard basis of diagonal matrices, then the positive roots
are $\omega_i$ for $1 \leq i \leq p$ and $\omega_j \pm \omega_i$ for 
$1 \leq i < j \leq p$, and in the case of $\su(p,q)$ and $\fsp(p,q)$, 
$2\omega_i$ for $1 \leq i \leq p$.

The root spaces $\n_{\om_j-\om_i}$, $\n_{\om_j+\om_i}$ and $\n_{2\om_k}$
are all contained in the space of matrices for which $C$, $E$ and $F$ are
all zero. The root spaces $\n_{\om_k}$, $1\leq k\leq p$, together span the space
$W$ of matrices for which $A$, $B$ and $D$ are zero and $C=E$.  For each $k$,
the root space $\n_{\om_k}$ is the set of such matrices for which all rows of $C$
except the $k$th one are zero.

Fix an integer $a$ with $1\leq a<m$, set $b=m-a$, and decompose $W$ into
subspaces $W=W_a+W_b$ where $W_a$, respectively $W_b$, is the space of
those matrices in $W$ for which the last $b$ columns, respectively first $a$
columns, of $C$ are zero.  We modify the Lie algebra $\s$ to construct a new Lie
algebra $\s_a$ by replacing all elements $X$ in $W_b$ by $\,X$ in
$\s^{\bC}$, leaving the spaces $W_a$ and the root spaces $2\omega_i$, $1 \leq i
\leq p$, and $\omega_j \pm \omega_i$, $1 \leq i < j \leq p$ unchanged.  The fact
that $\s_a$ is closed under bracket in $\s^{\bC}$ follows from the observations
that (i) $W$ commutes with all root spaces of the form $\n_{\om_j+\om_i}$
or $\n_{2\om_k}$, (ii) $W_b$ commutes with $W_a$, and (iii) $\a$ and the root
spaces of the form $\n_{\om_j-\om_i}$ normalize each of $W_a$ and $W_b$.

The simply-connected solvmanifold $S_a$ corresponding to $\s_a$ is
Einstein by Theorem \ref{thm1}.  $S_a$ has some 2-planes of positive
curvature and thus is not isometric to the symmetric space $S$.  Indeed, let $i$
and $j$ be any two positive integers with $i\leq a$ and $a <j \leq a+b$.  For
$1\leq k\leq m$, let $X_k$ be the element of $\n_{\om_k}$ for which $C$
has  the $(k,i)$ and $(k,j)$ entries equal to $\frac{1}{\sqrt{2}}$ and all
other entries zero.   Let $X'_k$ be the corresponding unit vectors in the
algebra $\s_a$.  Then for $k< l$, we have $[X'_k,X'_l]=0$.  Moreover,
$\langle U(X'_k,X'_k), U(X'_l,X'_l)\rangle=0$, but $U(X'_k,X'_l)$ is easily seen to
be a non-zero element of $\n_{\om_k-\om_l}$.  Thus by equation (\ref{sect}),
$K(X'_k,Y'_l)>0$.

Finally we verify that if $1 \leq a',a''<m$, then $S_a'$ and $S_a''$ are not
isometric unless either $a''=a'$ or $a''=m-a'$.  For simplicity, we will
just work with the case that $\bF=\R$, although similar arguments work in the
other cases.  In the original Lie algebra $\s$, consider two root spaces of the
form $\n_{\om_k}$, say $k=1$ and $k=2$.  Then the subspace $\h$ of
$\s$ spanned by $\n_{\om_1}$, $\n_{\om_2}$, and  $\n_{\om_1+\om_2}$ is
closed under brackets and is isomorphic to the $2m+1$-dimensional
Heisenberg algebra.  (See  Example \ref{heis}.)  Indeed, letting $X_i$ be the
element of $\n_{\om_1}$ for which the only non-zero entry of $C$ is the
$(1,i)$-entry, letting $Y_i$ be the analogous element of $\n_{\om_2}$, and letting
$Z$ be an appropriately chosen element of the one-dimensional space
$\n_{\om_1+\om_2}$, then the basis $\{X_1,\dots,X_m,Y_1,\dots,Y_m,Z\}$ of $\h$
satisfies the Heisenberg bracket relations $[X_i,Y_i]=Z$, with all other brackets
of basis vectors being zero.  In the modified algebra $\s_{a'}$, we replace
$X_i$ by $X'_i$ and $Y_i$ by $Y'_i$ with the new bracket relations
$[X'_i,Y'_i]=\pm Z$ where the sign is positive when $i\leq a'$ and negative when
$i> a'$.  The subspace $\h'$ of $\s_{a'}$ is still isomorphic to the
Heisenberg algebra, although the particular
basis $\{X'_1,\dots,X'_m,Y'_1,\dots,Y'_m,Z\}$ does not satisfy the canonical
Heisenberg bracket relations due to the introduction of negative signs.  Similar
statements hold in $\s_{a''}$.

The eigenspace $\n_{\om_2-\om_1}$ in $\s$ is one-dimensional and a choice
of basis vector $U$ satisfies $[U,X_i]=Y_i$ for all $i=1,\dots,m$; i.e.,
$\ad(U):\n_{\om_1}\to\n_{\om_2}$ is the canonical isomorphism $\sigma$ given by
$\sigma(X_i)=Y_i$. Now suppose that $\tau:\s_{a'}\to\s_{a''}$ is an isometry of
metric Lie algebras.  By Proposition \ref{isoroot}, $\tau$ induces an isomorphism
of the root system $\S$.  However, as one can see by a direct computation or by
referring to the Satake diagrams BI and DI on pages 532-533 of \cite{H}, the only
isomorphism of $\S$ is the identity.  It follows that $\tau$ restricts
to the identity on $\a$ and carries each root space $\n '_\al$ in $\n '$
to the corresponding root space $\n ''_\al$ in $\n ''$.  Applying this fact to
$\al=\om_1-\om_2$, we see that $\tau(U)=\pm U$ since $\tau$ is an isometry. 
Consequently, the restrictions $\tau:\n'_{\om_1}\to\n''_{\om_1}$ and
$\tau:\n'_{\om_2}\to\n''_{\om_2}$ satisfy $\tau\circ\sigma=\pm \sigma\circ\tau$. 
Finally, $\tau$ must restrict to an isomorphism between the Heisenberg algebras
$\h'$ and $\h''$.  It is easy to see that this is impossible unless $a''=a'$ or
$a''=m-a'$.

\subsection{SO(n,H)/U(n)}\label{SO(n,H)}
Our next non-compact symmetric space is dual to $\SO(2n)/\U(n)$, the compact
symmetric space of special orthogonal complex structures on $\R^{2n}$, with rank $m
= \lfloor\frac n2\rfloor$.  The dimension is $n^2 - n$.  On the Lie algebra level,
$$\so(n,\bH) = \left\{\pmatrix X & Y \\ -\overline Y & \overline X
\endpmatrix \in \gl(2n,\bC) \mid X = -X^t,~Y =
\overline Y^t \right\}.$$
Let $E_{ij}$ denote the skew-symmetric matrix with 1 in the $ij^{th}$ entry and
$-1$ in the $ji^{th}$ entry, and zeros elsewhere.  The Cartan subalgebra is
$$\a = \text{span}\left\{H_j = \tfrac{i}{\sqrt 2}(E_{2j-1,2j} -
E_{n+2j-1,n+2j})  \mid 1 \leq j \leq m\right\}.$$ 
Define $\omega_j \in \a^*$ by $\omega_j(H) = a_j$ for $H = i\sum_{j=1}^m a_j
H_j$, an arbitrary element of $\a$.  If $n$ is even  ($m = \frac n2$), the
positive roots of $\a$ are $\omega_j \pm \omega_k$ for $1 \leq j < k \leq m$, and
$2\omega_j$ for
$1 \leq j \leq m$.  If $n$ is odd  ($m = \frac{n-1}{2}$), we also have the
positive roots $\omega_j$ for $1 \leq j \leq m$.  We describe the root spaces
$\n_{\alpha}$:  
{\allowdisplaybreaks
\begin{align}
\n_{\omega_j \pm \omega_k} = \text{span}\{
A^{\pm}_{jk} = &\tfrac 12((E_{2j-1,2k-1} \mp E_{2j,2k} + E_{n+2j-1,n+2k-1}
\mp E_{n+2j,n+2k}) \notag
\\  +&i(\mp E_{2j-1,2k} - E_{2j,2k-1} \pm E_{n+2j-1,n+2k} +
E_{n+2j,n+2k-1})), \notag \\
B^{\pm}_{jk} = &\tfrac 12((E_{2j-1,2k} \pm E_{2j,2k-1} + E_{n+2j-1,n+2k} \pm
E_{n+2j,n+2k-1}) \notag
\\ +&i(\pm E_{2j-1,2k-1} - E_{2j,2k} \mp E_{n+2j-1,n+2k-1} +
E_{n+2j,n+2k})), \notag \\
C^{\pm}_{jk} = &\tfrac 12((E_{2j-1,n+2k} \mp E_{2j,n+2k-1} \mp E_{2k-1,n+2j}
+ E_{2k,n+2j-1}) \notag
\\ +&i(\mp E_{2j-1,n+2k-1} - E_{2j,n+2k} \pm E_{2k-1,n+2j-1} +
E_{2k,n+2j})), \notag \\
D^{\pm}_{jk} = &\tfrac 12((E_{2j-1,n+2k-1} + E_{2k-1,n+2j-1} \pm E_{2j,n+2k}
\pm E_{2k,n+2j}) \notag
\\ +&i(\pm E_{2j-1,n+2k} - E_{2j,n+2k-1} + E_{2k-1,n+2j} \mp
E_{2k,n+2j-1}))\}; \notag \\
\n_{2\omega_k} = \text{span}\{G_k = &\tfrac{1}{\sqrt 2}((E_{2k-1,n+2k-1} +
E_{2k,n+2k}) + i(E_{2k-1,n+2k} - E_{2k,n+2k-1}))\}; \notag \\
\text{(If $n$ odd)}\quad \n_{\omega_k} =  \text{span}\{X_k = &\tfrac{1}{\sqrt 2} 
((E_{2k,n} + E_{n+2k,2n}) +i(E_{2k-1,n} - E_{n+2k-1,2n})), \notag \\
Y_k = &\tfrac{1}{\sqrt 2}((E_{2k-1,n} + E_{n+2k-1,2n}) - i(E_{2k,n}
- E_{n+2k,2n})), \notag \\
Z_k = &\tfrac{1}{\sqrt 2}((E_{2k,2n} + E_{n,n+2k}) + i(E_{2k-1,2n}
-E_{n,n+2k-1})), \notag \\
W_k = &\tfrac{1}{\sqrt 2}((E_{2k-1,2n} + E_{n,n+2k-1}) - i(E_{2k,2n}
-E_{n,n+2k}))\}.\notag 
\end{align} }
These bases of the root spaces together with the basis above of $\a$ give an
orthonormal basis for the solvable subalgebra $\s$ of $\so(n,\bH)$.

We first consider the case of $n$ even.  The non-zero bracket relations $[X,Y]$
in $\n$ are given in the following table:

\begin{center}
\begin{tabular}{l||c|c|c|c|c|c|c|c|}
$Y\diagdown X$ & $A^+_{kl}$ & $A^-_{kl}$ & $B^+_{kl}$ & $B^-_{kl}$ &
$C^+_{kl}$ & $C^-_{kl}$ &
$D^+_{kl}$ & $D^-_{kl}$ \\
\hline\hline
$A^+_{kl}$ & 0 & 0 & 0 & 0 & 0 & 0 & 0 & $\sqrt 2\,G_k$ \\
\hline
$A^-_{jk}$ & $A^+_{jl}$ & $A^-_{jl}$ & $B^+_{jl}$ & $B^-_{jl}$ &
$C^+_{jl}$ & $C^-_{jl}$ & $D^+_{jl}$ & $D^-_{jl}$\\
\hline
$B^+_{kl}$ & 0 & 0 & 0 & 0 & 0 & $\sqrt 2\,G_k$ & 0 & 0 \\
\hline
$B^-_{jk}$ & $-B^+_{jl}$ & $B^-_{jl}$ & $A^+_{jl}$ & $-A^-_{jl}$ &
$-D^+_{jl}$ & $D^-_{jl}$ & $C^+_{jl}$ & $-C^-_{jl}$ \\
\hline
$C^+_{kl}$ & 0 & 0 & 0 & $-\sqrt 2\,G_k$ & 0 & 0 & 0 & 0 \\
\hline
$C^-_{jk}$ & $-C^+_{jl}$ & $C^-_{jl}$ & $D^+_{jl}$ & $-D^-_{jl}$ &
$A^+_{jl}$ & $-A^-_{jl}$ & $-B^+_{jl}$ & $B^-_{jl}$ \\
\hline
$D^+_{kl}$ & 0 & $-\sqrt 2\,G_k$ & 0 & 0 & 0 & 0 & 0 & 0 \\
\hline
$D^-_{jk}$ & $D^+_{jl}$ & $D^-_{jl}$ & $C^+_{jl}$ & $C^-_{jl}$ &
$-B^+_{jl}$ & $-B^-_{jl}$ & $-A^+_{jl}$ & $-A^-_{jl}$ \\
\hline
$G_l$ & 0 & $-D^+_{kl}$ & 0 & $-C^+_{kl}$ & 0 & $B^+_{kl}$ & 0 & $A^+_{kl}$ \\
\hline
\end{tabular}
\end{center}

We can modify $\s$ by replacing each of the basis elements 
$B^-_{jk}$, $C^-_{jk}$, $A^+_{jk}$, $D^+_{jk}$ and $G_k$ 
by their respective products with $\sqrt{-1}$, elements of $\s^{\bC}$, the
complexified Lie algebra.  With these changes, we obtain a new Lie algebra $\s'$ 
associated to $\s$, and the corresponding simply-connected Riemannian solvmanifold
$S'$ is Einstein.  $S'$ has some positive sectional curvature:  
Let $X =\frac{1}{\sqrt 2}((A^-_{j,j+1})' + (B^-_{j,j+1})')$ and let $Y =
\frac{1}{\sqrt 2} ((C^-_{j+1,j+2})' +(D^-_{j+1,j+2})')$, so that $[X,Y] = 0$ and
$U(X,Y) = 0$.  Then $U(X,X) = \frac{1}{\sqrt 2} (H_j - H_{j+1})$ and $U(Y,Y) =
\frac{1}{\sqrt 2} (H_{j+1} - H_{j+2})$. Hence $K(X,Y) = -\langle U(X,X),U(Y,Y)
\rangle >0$.  

Next consider the case of $n$ odd. In addition to the bracket relations above, we have
these additional non-zero Lie bracket relations $[X,Y]$:

\begin{center}
\begin{tabular}{l||c|c|c|c|c|c|c|c|}
$Y\diagdown X$ & $A^-_{kl}$ & $B^-_{kl}$ & $C^-_{kl}$ & $D^-_{kl}$
& $X_l$ & $Y_l$ & $Z_l$ & $W_l$ \\
\hline\hline
$X_l$ & $-X_k$ & $-Y_k$ & $-W_k$ & $Z_k$ & 0 & 0 & $\sqrt 2\,G_l$ & 0 \\
\hline
$Y_l$ & $-Y_k$ & $X_k$ & $-Z_k$ & $-W_k$ & 0 & 0 & 0 & $\sqrt 2\,G_l$ \\
\hline
$Z_l$ & $-Z_k$ & $-W_k$ & $Y_k$ & $X_k$ & $-\sqrt 2\,G_l$ & 0 & 0 & 0 \\
\hline
$W_l$ & $-W_k$ & $Z_k$ & $X_k$ & $Y_k$ & 0 & $-\sqrt 2\,G_l$ & 0 & 0 \\
\hline
$X_j$ & 0 & 0 & 0 & 0 & $A^+_{jl}$ & $-B^+_{jl}$ & $D^+_{jl}$ & $-C^+_{jl}$ \\
\hline
$Y_j$ &  0 & 0 & 0 & 0 & $-B^+_{jl}$ & $-A^+_{jl}$ & $C^+_{jl}$ & $D^+_{jl}$ \\
\hline
$Z_j$ & 0 & 0 & 0 & 0 & $-D^+_{jl}$ & $C^+_{jl}$ & $A^+_{jl}$ & $-B^+_{jl}$ \\
\hline
$W_j$ & 0 & 0 & 0 & 0 & $-C^+_{jl}$ & $-D^+_{jl}$ & $-B^+_{jl}$ &
$-B^+_{jl}$ \\
\hline
\end{tabular}
\end{center}

In this case, we modify $\s$ by replacing basis elements
$X_k$, $Z_k$, $B^+_{jk}$, $C^+_{jk}$, $B^-_{jk}$, $C^-_{jk}$
by their respective products with $\sqrt{-1}$. With this basis, we
obtain a new Lie algebra $\s'$ and thus an Einstein solvmanifold $S'$. Just as in
the case of $n$ even, $S'$ has positively curved sections.

\subsection{SL(n,H)/Sp(n)}\label{SL(n,H)}
Our final non-compact symmetric space has rank $n-1$; it is the dual space to
$\SU(2n)/\Sp(n)$, the compact symmetric space of special orthogonal quaternionic
structures on $\bC^{2n}$.  This space has dimension $2n^2-n-1$.  On the Lie
algebra level,
$$\sl(n,\bH) = \left\{ \pmatrix X & -\overline Y \\ Y & \overline X
\endpmatrix \mid X,Y \in \gl(n,\bC), \quad \tr(X + \overline X) = 0\right\}.$$
We have 
$$\a = \{diag(a_1,\dots,a_n,a_1,\dots,a_n) \in \gl(2n, \bC) \mid
\sum_{i=1}^n a_i =0\}.$$
Let $H = diag(a_1,\dots,a_n,a_1,\dots,a_n)$ be an arbitrary element of $\a$.  Then
define $\omega_j \in \a^*$ by $\omega_j(H) = a_j$.  The roots of $\s$ are
all of the form $\omega_j - \omega_k$, where $1 \leq j < k \leq n$.  Let $F_{ij}$
denote the matrix with a $1$ in the $ij^{th}$ entry and zeros elsewhere.  We
describe the root spaces of $\n$: 
{\allowdisplaybreaks
\begin{align}
\n_{\omega_j -\omega_k} = \{&A_{jk} = i\sqrt 2\,(F_{jk}-F_{n+j,n+k}),
~B_{jk} = i\sqrt 2\,(F_{j,n+k}+F_{n+j,k}), \notag \\
&C_{jk} = \sqrt 2\,(F_{j,n+k} - F_{n+j,k}), ~ D_{jk} = \sqrt 2\,(F_{jk} +
F_{n+j,n+k})\}. \notag
\end{align} }
We can represent the Lie brackets $[X,Y]$ in $\n$ by the following table
(notice that although we are shifting root spaces, one can think of the
multiplication of $A,B,C$ as quaternionic):
\begin{center}
\begin{tabular}{l||c|c|c|c|}
$Y\diagdown X$ & $A_{ij}$ & $B_{ij}$ & $C_{ij}$ & $D_{ij}$ \\
\hline\hline
$A_{jk}$ & $-\sqrt 2\,D_{ik}$ & $\sqrt 2\,C_{ik}$ & $-\sqrt 2\,B_{ik}$ &
$\sqrt 2\,A_{ik}$ \\
\hline
$B_{jk}$ & $-\sqrt 2\,C_{ik}$ & $-\sqrt 2\,D_{ik}$ & $\sqrt 2\,A_{ik}$ &
$\sqrt 2\,B_{ik}$ \\
\hline
$C_{jk}$ & $\sqrt 2\,B_{ik}$ & $-\sqrt 2\,A_{ik}$ & $-\sqrt 2\,D_{ik}$ &
$\sqrt 2\,C_{ik}$ \\
\hline
$D_{jk}$ & $\sqrt 2\,A_{ik}$ & $\sqrt 2\,B_{ik}$ & $\sqrt 2\,C_{ik}$ &
$\sqrt 2\,D_{ik}$ \\
\hline
\end{tabular}
\end{center}

We modify the Lie algebra in the following way: For all $1 \leq j < k \leq n$,
replace $A_{jk}$ by $\sqrt{-1}\,A_{jk}$ and $C_{jk}$ by $\sqrt{-1}\,C_{jk}$.  
This gives a new Lie algebra $\s'$.  The corresponding Einstein solvmanifold
$S'$  has some positive sectional curvature.   
Let $H_{ij} = X_{ii}-X_{jj}+ X_{n+i,n+i}-X_{n+j,n+j}$.  
Let $X =\frac{1}{\sqrt 2}(A_{ij}' + B_{ij}')$ and 
$Y = \frac{1}{\sqrt 2}(C_{jk}' + D_{jk}')$.Then
$[X,Y] = 0$, $U(X,Y) = 0$, and $U(X,X) = H_{ij}$, $U(Y,Y) = H_{jk}$ so that
$K(X,Y) = -\langle U(X,X),U(Y,Y) \rangle = -\langle H_{ij}, H_{jk}\rangle >0.$
%
%

\end{document}